\documentclass[12pt,leqno]{article}
\title{On prime factors of class number of  cyclotomic fields   }
\author{Roland Qu\^eme}
\usepackage{amsmath,amsthm,amstext,amsbsy}
\newtheorem{thm}{Theorem}[section]
\newtheorem{cor}[thm]{Corollary}

\newtheorem{lem}[thm]{Lemma}
\font\mathbb=msbm10

\newcommand{\Q}{\mbox{\mathbb Q}}
\newcommand{\Z}{\mbox{\mathbb Z}}

\newcommand{\modu}{\ \mbox{mod}\ }
\newcommand{\be}{\begin{equation}}
\newcommand{\ee}{\end{equation}}
\newcommand{\bd}{\begin{displaymath}}
\newcommand{\ed}{\end{displaymath}}
\newcommand{\bn}{\begin{enumerate}}
\newcommand{\en}{\end{enumerate}}
\date{2006 jan 09}
\begin{document}
\maketitle
\tableofcontents
\abstract Let $p$ be an odd prime. Let $K=\Q(\zeta)$ be the
$p$-cyclotomic number field. Let $v$ be a primitive root $\modu p$
and  $\sigma:\zeta\rightarrow\zeta^v$ be a $\Q$-isomorphism of the
extension $K/\Q$  generating the Galois group $G$ of $K/\Q$.
Following the conventions of Ribenboim in \cite{rib}, for $n\in\Z$,
the notation $v_n$ is understood  by $v_n=v^n\modu p$ with $1\leq
v_n\leq p-1$.
Let $P(X)=\sum_{i=0}^{p-2}X^iv_{-i}\in\Z[X]$
be the Stickelberger polynomial.
$P(\sigma)$ annihilates the class group $\mathbf C$ of $K$.
There exists a
polynomial $Q(X)\in\Z[G]$ such that $P(\sigma)\times(\sigma-v)=
p\times Q(\sigma)$ and such that  $Q(\sigma)$  annihilates the $p$-class
group $ C_p$ of $K$.
\clearpage
These  result allow:
\bn
\item
to describe  the structure of the relative class group $\mathbf C^-$,
\item
to give some explicit congruences in $\Z[v]\modu p$  for the $p$-class group of $K$ (the subgroup of exponent $p$ of $\mathbf C$),
\item
to give some explicit congruences in $\Z[v] \modu h$ for the $h$-class group
of $K$ for the prime divisors $h\not= p$ of the class number $h(K)$.
\item
We detail at the end the case of class number of quadratic and
biquadratic fields contained in the cyclotomic field $K$.
\item
As application, we give a MAPLE  algorithm which describes  the structure of the relative class group $\mathbf C^-$ of the cyclotomic field $K$
for
all the prime numbers $p<500$.
\en
This article is at elementary level.

{\it Remark:} we have not found in the literature some formulations
corresponding  to theorems
\ref{l603302} p. \pageref{l603302},
\ref{p607092} p. \pageref{p607092},
\ref{t611271} p. \pageref{t611271},
\ref{t607301} p. \pageref{t607301},
\ref{t607121} p. \pageref{t607121},
and
\ref{t607251} p. \pageref{t607251}
and to MAPLE algorithm described in section \ref{s701081} p \pageref{s701081}.
%
\section{Some definitions}\label{s601191}
In this section we give some definitions and  notations on cyclotomic
fields and  $p$-class group used in this paper.
\begin{enumerate}
\item
Let $p$ be an odd prime. Let ${\bf F}_p$ be the finite field of
cardinal $p$ and ${\bf F}_p^*$ its multiplicative group. Let $\zeta$
be a root of the polynomial equation $X^{p-1}+X^{p-2}+\dots+X+1=0$.
Let $K$ be the $p$-cyclotomic field $K=\Q(\zeta)$ and $O_K$ its
ring of integers. Let $K^+$ be the maximal totally real
subfield of $K$. Let $v$ be a primitive root $\modu p$. In this
paper, following Ribenboim conventions in  Ribenboim \cite{rib} for
any $n\in \Z$, we note $v_n=v^n\modu p$ with $1\leq v_n\leq p-1$.
Let $G$ be the Galois group of the extension $K/\Q$. Let $\sigma:
\zeta\rightarrow \zeta^v$  be a $\Q$-isomorphism of the extension
$K/\Q$  generating $G$. Let $\lambda=\zeta-1$. The prime ideal of
$K$ lying over $p$ is $\pi=\lambda O_K$.
\item
Let $\mathbf C$ be the class group of $K$.
Let $C_p$ be the $p$-class of $K$ (the  subgroup of exponent $p$ of $\mathbf C$ ).
Let $C_p^+$ be the $p$-class group of $K^+$. Then $C_p=C_p^+\oplus C_p^-$ where $C_p^-$ is called the relative $p$-class group.
Let $r^-$ be the rank of $C_p^-$.
\item
$C_p^-$ is the direct sum of $r^-$ subgroups  $\Gamma_k$  of order
$p$ annihilated by $\sigma-\mu_k\in {\bf F}_p[G]$ with $\mu_k\in{\bf
F}_p^*,\ \mu_k= v_{2m_k+1}$ where $m_k$ is a natural integer
$m_k,\quad 1\leq m_k \leq \frac{p-3}{2}$,
\be\label{e609231}
C_p^-=\oplus_{k=1}^{r^-} \Gamma_k.
\ee

\end{enumerate}
%
\section{On Kummer and Stickelberger relation}\label{s601192}
Stickelberger relation was already known by Kummer under the form of
Jacobi resolvents for the cyclotomic field $K$, see for instance
Ribenboim \cite{rib} (2.6) p. 119.
In this section we  derive some elementary properties  from  Stickelberger
relation.
\begin{enumerate}
\item
Let $q\not=p$ be an odd prime.
Let $\zeta_q$ be a root of the minimal polynomial equation $X^{q-1}+X^{q-2}+\dots+X+1=0$.
Let $K_q=\Q(\zeta_q)$ be the $q$-cyclotomic field.
Let $K_{pq}=\Q(\zeta_p,\zeta_q)$. Then $K_{pq}$ is the compositum $KK_q$.
The ring of integers of $K_{pq}$ is $O_{K_{pq}}$.
\item
Let $\mathbf q$ be a prime ideal of $O_K$ lying over the prime $q$.
Let $f$ be the order of $q\modu p$ and  $m=N_{K_p/\Q}(\mathbf q)=
q^f$. If $\psi(\alpha)=a$ is the image of $\alpha\in O_K$ under the
natural map $\psi: O_K\rightarrow O_K/\mathbf q$, then for
$\psi(\alpha)=a\not\equiv 0$ define a character $\chi_{\mathbf
q}^{(p)}$ on ${\bf F}_m=O_K/\mathbf q$ by
\begin{equation}
\chi_{\mathbf q}^{(p)}(a)={\{\frac{\alpha}{\mathbf q}\}}_p^{-1}=\overline{\{\frac{\alpha}{\mathbf q}\}}_p,
\end{equation}
where $\{\frac{\alpha}{\mathbf q}\}=\zeta^c$ for some natural integer $c$,
is the $p^{th}$ power residue character $\modu \mathbf q$.
We define the Gauss sum
\begin{equation}\label{e6012211}
g(\mathbf q)=\sum_{x\in{\bf F}_m}(\chi_{\mathbf q}^{(p)}(x)\times\zeta_q^{Tr_{{\bf F}_m/{\bf F}_q}(x)})\in O_{K_{pq}}.
\end{equation}
It follows that $\mathbf g(\mathbf q)\in O_{K_{pq}}$.
Moreover $g(\mathbf q)^p\in O_K$, see for instance Mollin \cite{mol} prop. 5.88 (c) p. 308.
\item
The Stickelberger relation is classically:
\begin{equation}\label{e512121}
\mathbf g(\mathbf q)^pO_K=\mathbf q^{S},
\end{equation}
with $S=\sum_{t=1}^{p-1} t\times \varpi_t^{-1}$,
where  $\varpi_t\in Gal(K/\Q)$ is given by $\varpi_t: \zeta\rightarrow \zeta^t$ (see for instance Mollin \cite{mol} thm. 5.109 p. 315).
\end{enumerate}
%

%
The four  following lemmas are derived in an elementary way from the Stickelberger relation.
\begin{lem}\label{l512151}
If $q\not\equiv 1\modu p$ then the Gauss sum $g(\mathbf q)\in \Z[\zeta]$.
\begin{proof}$ $
\begin{enumerate}
\item
Let $u$ be a primitive root $\modu q$. Let $\tau :\zeta_q\rightarrow
\zeta_q^u$ be a $\Q$-isomorphism generating $Gal(K_q/\Q)$. The
$\Q$-isomorphism $\tau$ is extended to a $K_p$-isomorphism of
$K_{pq}$ by $\tau:\zeta_q\rightarrow \zeta_q^u,\quad
\zeta_p\rightarrow \zeta_p$. Then  $g(\mathbf q)^p\in \Z[\zeta_p]$
and so
\begin{displaymath}
\tau(g(\mathbf q))^p=g(\mathbf q)^p,
\end{displaymath}
and it follows that there exists a natural integer $\rho$ with $\rho<p$ such that
\begin{displaymath}
\tau(g(\mathbf q))= \zeta_p^\rho\times  g(\mathbf q).
\end{displaymath}
Then $N_{K_{pq}/K}(\tau(g(\mathbf q)))=\zeta^{(q-1)\rho}\times N_{K_{pq}/K}(g(\mathbf q))$ and so  $\zeta^{\rho(q-1)}=1$.
\item
If $q\not\equiv 1\modu p$, it implies that $\zeta^\rho=1$ and so that $\tau(g(\mathbf q))=g(\mathbf q)$
and thus that $g(\mathbf q)\in O_K$.
\end{enumerate}
\end{proof}
\end{lem}
%
\begin{lem}\label{l12161}
Let $S=\sum_{t=1}^{p-1} \varpi_t^{-1}\times t$ where $\varpi_t$ is
the $\Q$-isomorphism of the extension $K/\Q$ given by
$\varpi_t:\zeta\rightarrow \zeta^t$ of $K$. Let
$P(\sigma)=\sum_{i=0}^{p-2} \sigma^i\times v_{-i}\in\Z[G]$. Then
$S=P(\sigma)$.
\begin{proof}
Let us consider one term $\varpi_t^{-1} \times t$. Then
$v_{-1}=v_{p-2}$ is a primitive root $\modu p$  and so there exists one and one $i$ such that
$t=v_{-i}$. Then $\varpi_{v_{-i}}:\zeta\rightarrow \zeta^{v_{-i}}$
and so $\varpi_{v_{-i}}^{-1}:\zeta\rightarrow\zeta^{v_i}$ and so
$\varpi_{v_{-i}}^{-1}=\sigma^i$ (observe that $\sigma^{p-1}\times
v_{-(p-1)}=1$), which achieves the proof.
\end{proof}
\end{lem}

%
\begin{lem}\label{l512171}
\begin{equation}
P(\sigma)=\sum_{i=0}^{p-2}\sigma^i\times v_{-i}=v_{-(p-2)}\times
\{\prod_{k=0,\ k\not=1}^{p-2}(\sigma-v_{k})\}+p\times R(\sigma),
\end{equation}
where $R(\sigma)\in\Z[G]$ with $deg(R(\sigma))<p-2$.
\begin{proof}
Let us consider the polynomial
$R_0(\sigma)=P(\sigma)-v_{-(p-2)}\times \{\prod_{k=0,\
k\not=1}^{p-2}(\sigma-v_{k})\}$ in ${\bf F}_p[G]$. Then
$R_0(\sigma)$ is of degree smaller than $p-2$ and  the two
polynomials $\sum_{i=0}^{p-2} \sigma^iv_{-i} $ and  $\prod_{k=0,\
k\not=1}^{p-2}(\sigma-v_{k})$ take a null value in ${\bf F}_p[G]$
when $\sigma$ takes  the $p-2$ different  values  $\sigma=v_k$ for
$k=0,\dots, p-2,\quad k\not= 1$. Then $R_0(\sigma)=0$ in ${\bf
F}_p[G]$ which leads to the result in $\Z[G]$.
\end{proof}
\end{lem}
%
\begin{lem}\label{l512165}
\begin{equation}\label{e512172}
P(\sigma)\times (\sigma-v)= p\times Q(\sigma),
\end{equation}
where $Q(\sigma)=\sum_{i=1}^{p-2}\delta_i\times \sigma^i\in\Z[G]$ is given by
\begin{equation}
\begin{split}
& \delta_{p-2}= \frac{v_{-(p-3)}-v_{-(p-2)}v}{p},\\
& \delta_{p-3}= \frac{v_{-(p-4)}-v_{-(p-3)}v}{p},\\
& \vdots\\
& \delta_i=\frac{v_{-(i-1)}-v_{-i} v}{p},\\
&\vdots\\
& \delta_1 = \frac{1-v_{-1}v}{p},\\
\end{split}
\end{equation}
with $-p<\delta_i\leq 0$.
\begin{proof}
We start of  the relation in $\Z[G]$
\begin{displaymath}
P(\sigma)\times(\sigma-v)= v_{-(p-2)}\times \prod_{k=0}^{p-2} (\sigma-v_k)+p\times R(\sigma)\times(\sigma-v)=p\times Q(\sigma),
\end{displaymath}
with $Q(\sigma)\in\Z[G]$ because  $\prod_{k=0}^{p-2} (\sigma-v_k)=0$ in ${\bf F}_p[G]$ and so
$\prod_{k=0}^{p-2} (\sigma-v_k)=p\times R_1(\sigma)$  in $\Z[G]$.
Then  we identify in $\Z[G]$ the  coefficients in the relation
\begin{displaymath}
\begin{split}
&(v_{-(p-2)}\sigma^{p-2}+v_{-(p-3)}\sigma^{p-3}+\dots+v_{-1}\sigma+1)\times(\sigma-v)=\\
&p\times (\delta_{p-2}\sigma^{p-2}+\delta_{p-3}\sigma^{p-3}+\dots+\delta_1\sigma+\delta_0),\\
\end{split}
\end{displaymath}
where $\sigma^{p-1}=1$.
\end{proof}
\end{lem}
\paragraph{Remarks:}
\begin{enumerate}
\item
Observe that we have more generally for the indeterminate $X$ the algebraic identity in $\Z[X]$
\be\label{e701031}
P(X)(X-v)= p\times Q(X)+v(X^{p-1}-1).
\ee
\item
Observe that, with our notations,  $\delta_i\in \Z,\quad i=1,\dots,p-2$, but generally $\delta_i\not\equiv 0\modu p$.
\item
We see also that $-p< \delta_i\leq 0$.
Observe also that $\delta_0=\frac{v^{-(p-2)}-v}{p}=0$.
\end{enumerate}
%
\section{Polynomial congruences $\modu p$  connected to the $p$-class group $C_p$}\label{s601193}
We give some  explicit polynomial congruences in $\Z[v]\modu p$ connected to  the relative  $p$-class group $C_p^-$  of $O_K$.
We apply successively the Stickelberger relation to prime ideals $\mathbf q$ of inertial degree $f=1$ and of inertial degree $f>1$.
We recall that $r^-$ is the $p$-rank of the relative $p$-class group $C_p^-$ of $K$.
%
\subsection{Stickelberger relation for prime ideals $\mathbf q$ of inertial degree $f=1$}
\begin{thm}\label{t512171}
 $\prod_{k=1}^{r^-}(\sigma-v_{2m_k+1})$ divides $Q(\sigma)$ in ${\bf F}_p[G]$ and
\begin{equation}\label{e512191}
Q(v_{2m_k+1})=\sum_{i=1}^{p-2} v_{(2m_k+1)\times i}\times(\frac{v_{-(i-1)}-v_{-i}\times v}{p})\equiv 0\modu p,
\end{equation}
where $m_k$ is defined in relation (\ref{e609231}) p. \pageref{e609231}.
\begin{proof}$ $
From Kummer, the group of ideal classes of $K$ is generated by the classes of prime ideals of degree $1$
(see for instance Ribenboim \cite{rib} (3A) p. 119).
Let $\mathbf q$ a prime ideal of inertial degree $1$ whose class $Cl(\mathbf q)\in C_p^-$  is annihilated by
$\sigma-\mu_k$ with $\mu_k=v_{2m_k+1}$.
We start of $g(\mathbf q)^p O_K= \mathbf q^{P(\sigma)}$ and so
$g(\mathbf q)^{p(\sigma-v)}O_K= \mathbf q^{P(\sigma)(\sigma-v)}=\mathbf q^{p Q(\sigma)}$ and thus
$g(\mathbf q)^{(\sigma-v)}O_K= \mathbf q^{ Q(\sigma)}$.
It can be shown that $g(\mathbf q)^{\sigma-v}\in K$, see for instance Ribenboim \cite{ri2} F. p. 440.
Therefore  $Q(\sigma)$ annihilates the ideal class  $Cl(\mathbf q)$ and the congruence follows.
\end{proof}
\end{thm}
\paragraph{Remarks:}
\begin{enumerate}
\item
Observe that $\delta_i$ can also be written in the form $\delta_i=-[\frac{v_{-i}\times v}{p}]$
where $[x]$ is the integer part of $x$, similar form also known in the literature.
\item
Observe that it is possible to get  other polynomials of $\Z[G]$ annihilating the relative $p$-class group $C_p^-$: for instance
from  Kummer's formula on Jacobi cyclotomic functions we induce other  polynomials $Q_d(\sigma)$ annihilating the
relative $p$-class group $C_p^-$  of $K$ : If $1\leq d\leq p-2$ define the set
\begin{displaymath}
I_d=\{i\ |\ 0\leq i\leq p-2, \quad v_{(p-1)/2-i}+v_{(p-1)/2-i+ind_v(d)}> p\}
\end{displaymath}
where the index $ind_v(d)$ is the minimal integer $s$ such that $d= v_s$.
Then the  polynomials $Q_d(\sigma)=\sum_{i\in I_d}\sigma^i$ for $d=1,\dots,p-2$ annihilate the $p$-class $C_p$ of $K$,
see for instance Ribenboim \cite{rib}  relations (2.4) and (2.5) p. 119.
\item
See also in a more general context Washington, \cite{was} corollary 10.15 p. 198.
\item
It is easy to verify the consistency of relation (\ref{e512191}) with the table of irregular primes and Bernoulli numbers in
Washington, \cite{was} p. 410.
\item
See section \ref{s701081} p. \pageref{s701081} with an algorithm improving this result.
\end{enumerate}
%
\subsection{Stickelberger relation for prime ideals $\mathbf q$ of inertial degree $f>1$}\label{s604192}
Let $\mathbf q$ be a prime ideal of $O_K$ with $Cl(\mathbf q)\in C_p$.
In this section we apply  Stickelberger relation  to the prime ideals $\mathbf q$ of  inertial degree $f>1$
with the method used for the prime ideals of inertial degree $1$ in section \ref{s601193} p. \pageref{s601193}.
Observe, from lemma \ref{l512151} p.
\pageref{l512151},  that  $f>1$ implies  that  $g(\mathbf q)\in O_K$, property used in this section
(by opposite $g(\mathbf q)\not\in K$ when $f=1$).

\paragraph{A definition:} we say that the prime ideal $\mathbf c$ of  a number field $M$ is $p$-principal if the component of the class group
$<Cl(\mathbf c)>$ in $p$-class group $D_p$ of $M$ is trivial.
%
\begin{thm}\label{l603302}
Let $q$ be an odd prime with $q\not=p$. Let $f$ be the order of $q\modu p$ and  $m=\frac{p-1}{f}$.
Let $\mathbf q$ be a  prime   ideal of $O_K$ lying over $q$ with $Cl(\mathbf q)\in C_p^-$.
If $f>1$ then
\bn
\item
$g(\mathbf q)\in O_K$  and $g(\mathbf q) O_K= \mathbf q^{P_1(\sigma)}$
where
\begin{equation}\label{e604013}
P_1(\sigma)=\sum_{i=0}^{m-1}(\frac{\sum_{j=0}^{f-1} v_{-(i+jm)}}{p})\times\sigma^i\in\Z[G].
\end{equation}
\item
There exists a natural integer $n$ with $1\leq n\leq p-2$ such that $\sigma-v_{n}$  divides $P_1(\sigma)$ in ${\bf F}_p[G]$
and $\sigma-v_{n}$ annihilates $Cl(\mathbf q)$.
\en
\begin{proof}$ $
\begin{enumerate}
\item

$N_{K/\Q}(\mathbf q)=q^f$
and
$\mathbf q=\mathbf q^{\sigma^m}=\dots=\mathbf q^{\sigma^{(f-1)m}}$.
From Stickelberger relation $g(\mathbf q)^p \Z[\zeta_p]=\mathbf q^{P(\sigma)}$
where $P(\sigma)=\sum_{i=0}^{m-1}\sum_{j=0}^{f-1}\sigma^{i+j m}v_{-(i+jm)}$.
Observe that, from hypothesis, $\mathbf q=\mathbf q^{\sigma^m}=\dots=\mathbf q^{\sigma^{(f-1)m}}$ so Stickelberger's relation implies that
$g(\mathbf q)^p O_K=\mathbf q^{P(\sigma)}$
with
\begin{displaymath}
P(\sigma)=\sum_{i=0}^{m-1}\sum_{j=0}^{f-1}\sigma^{i}v_{-(i+jm)}
= p\times \sum_{i=0}^{m-1} (\frac{\sum_{j=0}^{f-1} v_{-(i+jm)}}{p})\times  \sigma^i,
\end{displaymath}
where $(\sum_{i=0}^{f-1} v_{-(i+jm)})/p\in\Z$ because $v_{-m}-1\not\equiv 0\modu p$.
\item
Let $P_1(\sigma)=\frac{P(\sigma)}{p}$. Then from below $P_1(\sigma)\in \Z[G]$.
From lemma \ref{l512151} p.
\pageref{l512151} we know that  $f>1$ implies that   $g(\mathbf q)\in O_K$.
Therefore
\begin{displaymath}
g(\mathbf q)^p O_K=\mathbf q^{ p P_1(\sigma)}, \quad g(\mathbf q)\in O_K,
\end{displaymath}
and so
\begin{displaymath}
g(\mathbf q) O_K=\mathbf q^{P_1(\sigma)}, \quad g(\mathbf q)\in O_K.
\end{displaymath}
Then $P_1(\sigma)$ annihilates $Cl(\mathbf q)$ and  thus there exists $1 < n\leq p-2$  such that $(\sigma-v_n)\ |\ P_1(\sigma)$.
\end{enumerate}
\end{proof}
\end{thm}
%
\paragraph{Remarks}
\begin{enumerate}
\item
For $f=2$  the value of polynomial $P_1(\sigma)$ obtained from this lemma is $P_1(\sigma)=\sum_{i=1}^{(p-3)/2}\sigma^i$.
\item
Let $\mathbf q$ be a prime not principal ideal of inertial degree $f>1$ with $Cl(\mathbf q)\in C_p$. The {\bf two} polynomials of $\Z[G_p]$,
$Q(\sigma)=\sum_{i=0}^{p-2} (\frac{v_{-(i-1)}-v_{-i}v}{p})\times \sigma^i$ (see thm \ref{t512171})  and
$P_1(\sigma)=\sum_{i=0}^{m-1}(\frac{\sum_{j=0}^{f-1} v_{-(i+jm)}}{p})\times \sigma^i$   (see lemma \ref{l603302}) annihilate
the ideal class $Cl(\mathbf q)$.  When $f>1$ the lemma \ref{l603302} p. \pageref{l603302}
supplement the theorem \ref{t512171} p. \pageref{t512171}.
\item
This result explains that, when $f$ increases, the proportion of $p$-principal ideals $\mathbf q$ increases.
\end{enumerate}
%
It is possible to derive some explicit congruences in $\Z$ from this theorem.
\begin{cor}\label{l604012}
Let $q$ be an odd prime with $q\not=p$. Let $f$ be the order of $q\modu p$ and let $m=\frac{p-1}{f}$.
Let $\mathbf q$ be an  prime   ideal of $O_K$ lying over $q$.
Suppose that $f>1$.
\begin{enumerate}
\item
If  the ideal $\mathbf q$ is non $p$-principal   there exists  a natural integer $l,\ 1\leq l< m$ such that
\begin{equation}\label{e604015}
\sum_{i=0}^{m-1}(\frac{\sum_{j=0}^{f-1} v_{-(i+jm)}}{p})\times v^{lfi}\equiv 0\modu p,
\end{equation}
\item
If for all natural integers $l$  such that $1\leq l<m$
\begin{equation}\label{e604011}
\sum_{i=0}^{m-1}(\frac{\sum_{j=0}^{f-1} v_{-(i+jm)}}{p})\times v^{lfi}\not\equiv 0\modu p,
\end{equation}
then $\mathbf q$ is $p$-principal
\end{enumerate}
\begin{proof}$ $
\begin{enumerate}
\item
Suppose that $\mathbf q$ is not $p$-principal.
Observe at first that congruence (\ref{e604015}) with $l=m$ should imply that
$\sum_{i=0}^{m-1}(\sum_{j=0}^{f-1} v_{-(i+jm)})/p)\equiv 0\modu p\ $ or
\bd
\sum_{i=0}^{m-1}\sum_{j=0}^{f-1} v_{-(i+jm)}\equiv 0\modu p^2,
\ed
 which is not possible because
$v_{-(i+jm)}=v_{-(i^\prime+ j^\prime m)}$ implies that $j=j^\prime$ and $i=i^\prime$ and so that
$\sum_{i=0}^{m-1}\sum_{j=0}^{f-1} v_{-(i+jm)}=\frac{p(p-1)}{2}$.
\item
The polynomial $P_1(\sigma)$ of lemma \ref{l603302} annihilates the
non $p$-principal ideal $\mathbf q$
in ${\bf F}_p[G]$  only if there exists $\sigma-v_n$ dividing $P_1(\sigma)$ in ${\bf F}_p[G]$. From $\mathbf q^{\sigma^m-1}=1$ it follows also that
$\sigma-v_n\ |\ \sigma^m-1$. But $\sigma-v_n\ |\ \sigma^m-v_{nm}$ and so $\sigma-v_n\ |\ v_{nm}-1$, thus $nm\equiv 0\modu p-1$,  so
$n\equiv 0 \modu f$ and $n=lf$ for some $l$.
Therefore if $\mathbf q$ is non $p$-principal there exists  a natural integer $l,\ 1\leq l< m$ such that
\begin{equation}\label{e604014}
\sum_{i=0}^{m-1}(\frac{\sum_{j=0}^{f-1} v_{-(i+jm)}}{p})\times v^{lfi}\equiv 0\modu p,
\end{equation}
\item
The relation (\ref{e604011}) is an immediate consequence of previous part of the proof.
\end{enumerate}
\end{proof}
\end{cor}
%
%
\subsection{Polynomial congruences $\modu p^2$ connected to the $p^2$-class group $C_{p^2}$}
Let $C_{p^2}^-$ be the subgroup of exponent $p^2$ (so with elements of order dividing $p^2$) of the relative class group
$\mathbf C^-$ of $K$.
\bn
\item
We have seen in relation (\ref{e609231}) p. \pageref{e609231} that
the relative $p$-class group $C_p^-$ can be seen as a direct sum $C_p^-=\oplus_{i=1}^{r^-} \Gamma_i$
where $\Gamma_i$ is a cyclic subgroup of $C_p^-$ annihilated by $\sigma-\mu_i$ with $2\leq \mu_i\leq p-2$.
\item
$C_{p^2}^-$ can be seen as a direct sum
\be\label{e612061}
C_{p^2}^-=\oplus_{i=1}^{r^-} \Delta_i,
\ee
where $\Delta_i$ is a cyclic group with $\Gamma_i\subset\Delta_i$ and whose order divides $p^2$  .
\item
Suppose that $\Delta_i$ is of order $p^2$.
Show that  $\Delta_i$  is  annihilated by $\sigma-(\mu_i+a_ip)$ with  $a_i$  natural integer $1\leq a_i\leq p-1$:
\bn
\item
From Kummer,  there exist some  prime ideals $\mathbf Q_i$ of $O_K$ with $Cl(\mathbf Q_i)\in\Delta_i$, $Cl(\mathbf Q_i^p)\in \Gamma_i$
and $<\mathbf Q_i^{p(\sigma-\mu_i)}>$  principal as seen in  previous sections.
Therefore $<\mathbf Q_i^{\sigma-\mu_i}>$
is of order $p$.
\item
 $Cl(\mathbf Q_i^{\mu_i})\ \in \Delta_i$. In the other hand  $Cl(\mathbf Q_i^\sigma)\in \Delta_i$:
if not
$Cl(\mathbf Q_i^\sigma)$ should have at least one  component $\mathbf c_j\in \Delta_j,\ j\not=i$ and so
$Cl(\mathbf Q_i^{p\sigma})$ should have a component $\mathbf c_j^p\in \Delta_j^p=\Gamma_j,\ j\not=i$, contradiction.
Therefore $Cl(\mathbf Q_i^{\sigma-\mu_i})\in \Delta_i$.
Then, from $\mathbf Q_i^{p(\sigma-\mu_i)}$ principal,
it follows that $Cl(\mathbf Q_i^{\sigma-\mu_i})\in \Gamma_i$ because $\Delta_i$ is cyclic of order $p^2$.
\item
Thus there exists $a_i,\ 1\leq a_i\leq p-1,$ such that $Cl(\mathbf Q_i^{\sigma-\mu_i})= Cl(\mathbf Q_i^{a_i p})$
and so
$\mathbf Q_i^{\sigma-\mu_i-a_i}$ is principal.
\en
\en
In this section we examine  the case of subgroups $\Delta_i$ of order $p^2$.
Let us note $\Delta$ for one of this groups annihilated by $\sigma-(\mu+ap),\ a\not=0$.
%
\begin{thm}\label{t612061}
$\mu$ verifies the two congruences
\be\label{e612062}
\begin{split}
&\sum_{i=0}^{p-2}\mu^{i}\delta_i\equiv 0 \modu p \ \mbox{with\ }\delta_i=\frac{v_{-(i-1)}-v_{-i}v}{p},\\
&\sum_{i=0}^{p-2}\mu^{p-2+i}\delta_i+(\mu^{p-1}-1)\sum_{i=1}^{p-2} i\mu^{i-1}\delta_i\equiv 0\modu p^2.\\
\end{split}
\ee
\begin{proof} $ $
\bn
\item
There exists  prime ideals  $\mathbf Q$ of $O_K$ with $Cl(\mathbf Q)\in \Delta$, hence  $\mathbf Q^{p^2}$   principal
 and $\mathbf Q^p$ not principal.
 From Stickelberger relation
$ g(\mathbf Q)^pO_K=\mathbf Q^{P(\sigma)}$
where $P(\sigma)$ has been defined in lemma \ref{l12161} p. \pageref{l12161}.
Then
$ g(\mathbf Q)^{p(\sigma-v)}O_K=\mathbf Q^{P(\sigma)(\sigma-v)}$,
hence from lemma \ref{l512165} p. \pageref{l512165} we get
$ g(\mathbf Q)^{p(\sigma-v)}O_K=\mathbf Q^{p Q(\sigma)}$,
hence
$ g(\mathbf Q)^{(\sigma-v)}O_K=\mathbf Q^{ Q(\sigma)}$.
We know, for instance  from Ribenboim \cite{ri2} F. p. 440  that $g(\mathbf Q)^{\sigma-v}\in K$
so
$\mathbf Q^{Q(\sigma)}$ is principal.
But $\mathbf Q^{\sigma-(\mu+ap)}$ is principal  hence
$\mathbf Q^{Q(\mu+ap)}$ is principal,
and thus
\be\label{e612071}
Q(\mu+ap)\equiv 0\modu p^2.
\ee
\item
From lemma \ref{l512165} p. \pageref{l512165},
$Q(\sigma)=\sum_{i=0}^{p-2}\sigma^i\delta_i$ where
$\delta_i=\frac{v_{-(i-1)}-v_{-i}v}{p}$.
From relation (\ref{e612071})
\bd
\sum_{i=0}^{p-2}(\mu+ap)^i\delta_i\equiv 0\modu p^2,
\ed
hence
\bd
\sum_{i=0}^{p-2}\mu^i\delta_i+ap\sum_{i=1}^{p-1} i\mu^{i-1}\delta_i\equiv 0\modu p^2.
\ed
From theorem \ref{t512171} p. \pageref{t512171} applied to the ideal $\mathbf Q^p\in\Gamma$ of order $p$,
\bd
\sum_{i=0}^{p-2}\mu^i\delta_i\equiv 0\modu p.
\ed
$\mathbf Q^{\sigma^{p-1}-1}$ is principal,
therefore $\mathbf Q^{(\mu+ap)^{p-1}}-1$ is principal and so
$(\mu+ap)^{p-1}-1\equiv 0\modu p^2$, hence
$\mu^{p-1}+(p-1)\mu^{p-2}ap-1\equiv 0\modu p^2$,
hence
$\mu^{p-1}-1-\mu^{p-2}ap\equiv 0\modu p^2$,
hence
$a\equiv \frac{\mu^{p-1}-1}{\mu^{p-2}p}\modu p$,
and so
\bd
\sum_{i=0}^{p-2}\mu^i\delta_i+\frac{\mu^{p-1}-1}{\mu^{p-2}}\sum_{i=1}^{p-1} i\mu^{i-1}\delta_i\equiv 0\modu p^2
\ed
and finally
\be\label{e612073}
\sum_{i=0}^{p-2}\mu^{p-2+i}\delta_i+(\mu^{p-1}-1)\sum_{i=1}^{p-2} i\mu^{i-1}\delta_i\equiv 0\modu p^2,
\ee
which achieves the proof.
\en
\end{proof}
\end{thm}
%
\paragraph{Example:}
\bn
\item
This congruence $\modu p^2$ is  valid for  no irregular prime numbers $p<4001$ with rank $r$ of $C_p$ verifying $r>1$ (verified with a MAPLE program).
Therefore the class group of $K$ has no  cyclic subgroups  of order $p^2$ for the primes

$p=157, 353, 379, 467, 491, 547, 587, 617, 647, 673, 691, 809, 929, 1151, 1217,$

$1291, 1297, 1307, 1663,1669,
1733, 1789, 1847, 1933, 1997, 2003, 2087,$

$ 2273, 2309, 2371, 2383, 2423, 2441, 2591, 2671, 2789, 2909, 2939, 2957,$

$ 3391, 3407,
3511, 3517, 3533, 3539, 3559, 3593, 3617, 3637, 3833, 3851, 3881$.
\item
(see table of irregular prime in Washington \cite{was} p. 410).
Our result is consistent for the primes  $p=157, 353, 467, 491$ with Schoof \cite{sch}  table 4.2 p. 1239 describing structure of
class groups of some
$p$-cyclotomic fields.
\en
%
%
\section{On prime factors $h\not=p$ of the class number  of the $p$-cyclotomic field}\label{s607121}
In previous sections we considered the  relative $p$-class group $C_p^-$ of $K$.
By opposite, in this section we apply Stickelberger  relation to all the  primes $h\not=p$ dividing the class number $h(K)$.
A first subsection is devoted to the general case of the relative class group $\mathbf C^-$ of $K$, a second to
the class group of the quadratic subfield of $K$
and a third subsection to the class group to the biquadratic subfield of $K$ when $p\equiv 1\modu 4$.
%
\subsection{The general case}
\bn
\item
The class group $\mathbf C$ of $K$ is the direct sum of the class group $\mathbf C^+$ of the maximal totally real subfield $K^+$ of $K$ and
of the relative class group $\mathbf C^-$ of $K$.
\item
Remind that $v$ is a primitive root $\modu p$ and that $v_n$ is to be be understood as $v^n\modu p$ with $1\leq v_n\leq p-1$.
Let $h(K)$ be the class number of $K$. Let $h\not= p$ be an odd   prime dividing $h(K_p)$, with $v_h(h(K))=\beta$.
Let $d=Gcd(h-1,p-1)$.
Let $\mathbf C(h)$ be the $h$-Sylow subgroup of the class group of $K$ of order $h^\beta$.
Then
$\mathbf C(h)=\oplus_{j=1}^{\rho}\mathbf  C_{j}(h)$ where $\rho$ is the $h$-rank of the abelian group $\mathbf  C(h)$
of order $h^\beta$ and $\mathbf C_j(h)$
are cyclic groups of order $h^{\beta_j}$ where $\beta=\sum_{i=1}^\rho\beta_j$.
\item
From Kummer (see for instance Ribenboim \cite{rib} (3A) p. 119), the prime ideals of $O_K$ of inertial degree $1$
generate the ideal class group. Therefore there exist in the subgroup $C_h$ of exponent $h$ of $\mathbf C(h)$
some prime ideals $\mathbf q$   of inertial degree $1$  such that
$Cl(\mathbf q)\in \oplus_{j=1}^\rho c_j$ where $c_j$ is a cyclic group of order $h$
and $Cl(\mathbf q)\not \in \oplus_{j\in J}^\rho c_j$ where $J$ is a
strict subset of $\{1,2,\dots,\rho\}$.
\item
Let $P(\sigma)=\sum_{k=0}^{p-2}\sigma^k v_{-k}$ be the Stickelberger polynomial.
From lemma \ref{l12161} p. \pageref{l12161} Stickelberger relation is
$\mathbf q^{P(\sigma)}=g(\mathbf q)^p O_K$ where $g(\mathbf q)^p\in O_K$.
Therefore $\mathbf q^{P(\sigma)}$ is  principal, a fortiori is $\mathbf C(h)$-principal (or $Cl(\mathbf q)^{P(\sigma)}$
has a trivial component in $\mathbf C(h)$).
There exists a minimal polynomial $V(X)\in {\bf F}_h(X)$ of degree $\delta\leq \rho$ such that
$\mathbf q^{V(\sigma)}$  is $\mathbf C(h)$-principal, if not the remainder of the division of $P(X)$ by $V(X)$ of degree smaller than $\delta$
would annihilate also $\mathbf q$.
Therefore the irreducible polynomial $V(X)$ divides $P(X)$ in ${\bf F}_h[X]$ for the indeterminate $X$.
\item
If $Cl(\mathbf q)\in \mathbf C^-$ then $\mathbf q^{\sigma^{(p-1)/2}+1}$ is principal.
\item
Let $D(X)\in{\bf F}_h(X)$ defined by $D(X)=Gcd(P(X), X^{(p-1)/2}+1)\modu h$.
\en
%
We  obtain  the following:
%
\begin{thm}\label{p607092}
Suppose that the prime $h\not=p$ divides the class number $h(K)$.
Let $D(X)= Gcd (P(X), X^{(p-1)/2}+1)\modu h$.
Then
\bn
\item
$V(X)$ divides $D(X)$ in ${\bf F}_h[X]$.
The $h$-rank  $\rho$ of $\mathbf C(h)$ is greater or equal to  the degree $\delta$ of $V(X)$ and  $h^\rho\ |\ h(K)$.
\item
If $h$ is coprime with $p-1$ and   with the class number $h(E)$ of all intermediate fields $\Q\subset E\subset K,\ E\not= K$
then
$f\ |\ \rho$ where $f$ is the order of $h\modu p-1$.
\en
\begin{proof}$ $
\bn
\item
Reformulation of previous paragraph.
\item
Immediate consequence of theorem 10.8 p. 187 in Washington \cite{was} for the cyclic extension $K/E$.
\en
\end{proof}
\end{thm}
\paragraph{Remark:}
See the section \ref{s701081} p. \pageref{s701081} for a MAPLE program applying  the  theorem \ref{p607092} p. \pageref{p607092}.
%
%
\begin{cor}\label{t607072}
If $\mathbf C(h)$ is cyclic  then:
\bn
\item
$V(X)=X-\nu$ with $\nu\in{\bf F}_h^*$.
\item\label{i609242}
In ${\bf F}_h(X)$
\be\label{e607101}
\begin{split}
&V(X)\ |\  X^d-1,\ d=gcd(h-1,p-1),\\
& V(X)\ |\ \{\sum_{i=0}^{d-1}  X^i\times\frac{\sum_{j=0}^{(p-1)/d-1} v^{-(i+jd)}}{p}\}.\\
\end{split}
\ee
\item\label{l611261}
Let $M$ be the smallest subfield of $K$ such that $h\ |\ h(M)$.
Let $n=[M:\Q]$.
If $h$ is coprime with $n$ then   $h-1\equiv 0 \modu n$.
\en
\begin{proof}$ $
\bn
\item
$\rho=1$ implies that $V(\sigma)=\sigma-\nu$.
\item
$X-\nu\ |\ X^{h-1}-\nu^{h-1}$  and so $\sigma^{h-1}-\nu^{h-1}$ annihilates $\mathbf C(h)$.
From $\nu^{h-1}\equiv 1\modu h$ it follows that $\sigma^{h-1}-1$ annihilates $\mathbf C(h)$.
$\sigma^{p-1}-1$ annihilates $\mathbf C(h)$  and so
$\sigma^d-1$ annihilates $\mathbf C(h)$ and  so $X-\nu$ divides $X^d-1$ in ${\bf F}_h[X]$.
Then apply theorem  \ref{p607092}.
Observe that $\frac{\sum_{j=0}^{(p-1)/d-1} v^{-(i+jd)}}{p}\equiv 0\modu p$ and that we have assumed that the prime $h\not=p$ in this section.
\item
Let $E$ be the smallest  intermediate field $E,\ \Q\subset E\subset K$ with $h\ |\ h(E)$ and $[E:\Q]=n$.
Then we apply theorem 10.8 p. 187 of Washington \cite{was} to the cyclic extension $M/\Q$, thus
$h\equiv 1 \modu n$ and  $n\ |\ p-1$ because $f\ |\ \rho$ where $\rho=1$.
\en
\end{proof}
\end{cor}
%
It is possible to enlarge  previous results  with another annihilation polynomial:
\begin{thm}\label{t611271}
The polynomial $\Pi(\sigma)=\sum_{i=0,\ v_{-i}\ \mbox{even}}^{p-2}\sigma^i$
annihilates the non-$p$-part of the class group  $\mathbf C$ of the
$p$-cyclotomic field $K$.
\begin{proof}$ $

\bn
\item
Apply Stickelberger relation to field $\Q(\zeta_{2p})=\Q(\zeta_p)$.
Let $\varpi_{2t+1}: \zeta_{2p}\rightarrow\zeta_{2p}^{2t+1}$.
The Stickelberger polynomial can be written
\bd
S_2=\sum_{2t+1=1,\ t\not=p}^{2p-1}\varpi_{2t+1}^{-1}\times(2t+1).
\ed
\item
Observe at first that $\zeta_{2p}=-\zeta_p$.

\underline{If $t>(p-1)/2$}

then $\varpi_{2t+1}: \zeta_{2p}\rightarrow -\zeta_{2p}^{2t+1-p}$,
hence
$\varpi_{2t+1}: -\zeta_p\rightarrow -(-\zeta_p)^{2t+1-p}$,
hence
$\varpi_{2t+1}: -\zeta_p\rightarrow -(\zeta_p)^{2t+1-p}$ because $2t+1-p$ is even,
hence
$\varpi_{2t+1}: \zeta_p\rightarrow \zeta_p^{2t+1-p}$,
hence
$\varpi_{2t+1}=\varpi_{2t+1-p}$,
hence
\bd
\varpi_{2t+1}^{-1}\times (2t+1)=\varpi_{2t+1-p}^{-1}\times (2t+1-p)+p\times\varpi_{2t+1-p}^{-1}.
\ed
\item
The Stickelberger polynomial is
\bd
S_2=\sum_{t=0}^{(p-3)/2}\varpi_{2t+1}^{-1}\times(2t+1)
+\sum_{t=(p+1)/2}^{p-1}\varpi_{2t+1}^{-1}\times(2t+1),
\ed
hence
\bd
S_2=\sum_{t=0}^{(p-3)/2}\varpi_{2t+1}^{-1}\times(2t+1)
+\sum_{t=(p+1)/2}^{p-1}\varpi_{2t+1-p}^{-1}\times(2t+1-p)
+p\sum_{t=(p-1)/2}^{p-1}\varpi_{2t+1-p},
\ed
hence
\bd
S_2=\sum_{t=1}^{p-1}\varpi_t^{-1} t+p\times\sum_{t=(p-1)/2}^{p-1}\varpi_{2t+1-p},
\ed
hence, with $P(\sigma)$ defined in lemma \ref{l12161} p. \pageref{l12161},
\bd
S_2=P(\sigma)+p \times\sum_{t=(p-1)/2}^{p-1} \varpi_{2t+1-p},
\ed
With $2t+1-p=v_{-i}$  we get
\be\label{e611271}
S_2=P(\sigma)+p\times \sum_{i=0,\  v_{-i} \mbox{\ even}}^{p-2}\sigma^{i}.
\ee
\item
The polynomial $P(\sigma)$ annihilates the class group $\mathbf C$ of $K$.
Therefore the polynomial $p \times\sum_{i=0,\ v_{-i}\ \mbox{even}}^{p-2}\sigma^i$ annihilates also $\mathbf C$.
If $h\not=p$ then $\Pi(\sigma)$ annihilates $\mathbf C(h)$, which achieves the proof.
\en
\end{proof}
\end{thm}
%
\paragraph{Remark:}

Numerical MAPLE computations seem to show more :
the polynomial
\bd
Gcd(\Pi(\sigma),\sigma^{(p-1)/2}+1)=\frac{\sigma^{(p-1)/2}+1}{\sigma-v}=\prod_{m=1}^{(p-3)/2}(\sigma-v_{2m+1}).
\ed
Therefore $\Pi(\sigma)$ annihilates also the relative $p$-class group  $C_p^-$.
%

%
\begin{lem}\label{l607262}
Let $E$ be a subfield of $K$ with $[K:\Q]=d$.
Let $h$ be an odd  prime number dividing $h(L)$.
Then in ${\bf F}_h[X]$
\be\label{e607263}
\begin{split}
& V(X)\ |\  \sum_{i=0}^{d-1} X^i,\\
& V(X)\ |\ \{\sum_{i=0}^{d-1}  X^i\times\frac{\sum_{j=0}^{(p-1)/d-1} v^{-(i+jd)}}{p}\}.\\
\end{split}
\ee
\begin{proof}
$\sigma^d-1$ annihilates $\mathbf C(h)$.
The Stickelberger polynomial
\bd
P(X)=\sum_{i=0}^{p-2}X^i v_{-i}=p\times\sum_{i=0}^{d-1}  X^i\times\frac{\sum_{j=0}^{(p-1)/d-1} v^{-(i+jd)}}{p},
\ed
and from $p\not=h$ it follows that
\bd
V(X)\ |\ \{\sum_{i=0}^{d-1}  X^i\times\frac{\sum_{k=0}^{(p-1)/d-1} v^{-(i+jd)}}{p}\}
\ed
in ${\bf F}_p[G]$.
\end{proof}
\end{lem}
%
\paragraph{Remark:} Compare lemma \ref{l607262} for $h\not=p$ with lemma  \ref{l603302} p. \pageref{l603302} proved when $h=p$.
%

%
\subsection{The case of complex quadratic fields contained in $K$}

In this paragraph  we formulate directly  previous result when $h$ divides the class number of the complex quadratic field
$\Q(\sqrt{-p})\subset K, \ p\equiv 3 \modu 4,\ p\not=3$.
%
\begin{thm}\label{t607121} {Hilbert 145 theorem}

Suppose that $p\equiv 3\modu 4,\ p\not= 3$.
If $h$ is an odd  prime with $h\ |\ h(\Q(\sqrt{-p}))$  then
\begin{equation}\label{e607091}
\sum_{i=0}^{p-2}  (-1)^i v_{-i}\equiv 0\modu h.
\end{equation}
\begin{proof}
Let $\mathbf Q$ be the prime of $\Q(\sqrt{-p})$ lying above $\mathbf q$. The ideals $\mathbf Q\not=\sigma(\mathbf Q)$
and so $\mathbf Q^{\sigma+1}$ is principal because $\mathbf Q^{\sigma^2}=\mathbf Q$.
Therefore  $\mathbf Q^{\sum_{i=0}^{p-2} (-1)^i v_{-i}}$ is principal and
$\sum_{i=0}^{p-2}(-1)^i v_{-i}\equiv 0\modu h$.
\end{proof}
\end{thm}
%
\paragraph{Remarks:}
\bn
\item
The theorem \ref{t607121}  can also be obtained from  Hilbert Theorem 145 see Hilbert \cite{hil} p. 119.
See also Mollin, \cite{mol}  theorem 5.119 p. 318.
\item
From lemmas \ref{l512165} p. \pageref{l512165} we could prove similarly:

Suppose that $p\equiv 3\modu 4,\ p\not= 3$.
If $h$ is an odd prime with $h\ |\ h(\Q(\sqrt{-p}))$  then
\begin{equation}\label{e607091}
2\times\sum_{i=0}^{(p-3)/2}  (-1)^i v_{-i}-p\equiv 0\modu h.
\end{equation}
\item
Numerical evidences easily computable   show more:
If $p\not=3$ is prime with $p\equiv 3\modu 4$ then the class number $h(\Q(\sqrt{-p})$ verifies
\be
h(\Q(\sqrt{-p})=-\frac{\sum_{i=0}^{p-2} (-1)^i v_{-i}}{p}.
\ee
This result has been proved by Dirichlet by analytical number theory, see Mollin  remark 5.124 p. 321.
It is easy to verify  this formula, for instance in tables of  class numbers of complex quadratic fields in :
\bn
\item
H. Cohen \cite{coh} p. 502- 505, all the table for $p\leq 503$.
\item
in Wolfram table of quadratic class numbers  \cite{wol} for large $p$.
\en
\item
When $p\equiv 1 \modu 4$ this method cannot be applied to the quadratic field $\Q(\sqrt{p})\subset K$ because
$\sum_{i=0}^{p-2}(-1)^i v_{-i}$ is trivially null.
\en
%
\begin{thm}\label{t607301}
Suppose that $p\equiv 3\modu 4,\ p\not= 3$.
If $h$ is an odd  prime with $h\ |\ h(\Q(\sqrt{-p}))$  then
\begin{equation}\label{e607091}
\begin{split}
&\sum_{i=0,\ v_{-i} \mbox{\ even}}^{p-2}  (-1)^i \not= 0,\\
&\sum_{i=0,\ v_{-i} \mbox{\ even}}^{p-2}  (-1)^i \equiv 0\modu h.\\
\end{split}
\end{equation}
\begin{proof}$ $
We apply previous theorem \ref{t611271} p. \pageref{t611271} observing that in that case $\sigma^2(\mathbf Q)=\mathbf Q$ when $\mathbf Q$ is a non-principal ideal of
the quadratic field $\Q(\sqrt{-p})$, hence $\sigma+1$ annihilates $\mathbf C(h)$
and $\sum_{i=0,\ v_{-i} \mbox{\ even}}^{p-2}  (-1)^i \equiv 0\modu h$.
This sum has $\frac{p-1}{2}$ elements, thus of odd cardinal  and cannot be null.
\end{proof}
\end{thm}
%
\begin{thm}\label{t607121}
Suppose that $p\equiv 3\modu 4,\ p\not= 3$.
Let $\delta$ be an integer $1\leq \delta\leq p-2$.
Let $I_\delta$ be the set
\be\label{e607211}
\begin{split}
& I_\delta=\{i\ |\ 0\leq i\leq p-2,\  v_{(p-1)/2-i}+v_{(p-1)/2-i+ind_v(\delta)}>p\},
\end{split}
\ee
where, as seen above,  $ind_v(\delta)$ is the notation index of $\delta$ relative to $v$.
If $h$ is an odd  prime with $h\ |\ h(\Q(\sqrt{-p}))$  then
\begin{equation}\label{e607091}
\begin{split}
& \sum_{i\in I_\delta}  (-1)^i\not= 0,\\
& \sum_{i\in I_\delta}  (-1)^i\equiv 0\modu h.\\
\end{split}
\end{equation}
\begin{proof}
$I_\delta$ has an odd cardinal.
Then see relation (\ref{e607211}).
\end{proof}
\end{thm}
%
\paragraph{Remark:}
\bn
\item
Observe that results of theorems \ref{t607301} and \ref{t607121} are consistent with existing tables
of quadratic fields, for instance Arno, Robinson, Wheeler \cite{arn}.
Numerical verifications seem to show more :
\be\label{e67301}
\sum_{i=0,\ v_{-i} \mbox\ even}^{p-2} (-1)^i\equiv 0 \modu h(\Q(\sqrt{-p}).
\ee
\item
Observe that if $p\equiv 1\modu 4$ then $\sum_{i=0, v_{-i}\mbox{\ even}}^{p-2}(-1)^i=0$.
\en
%
\subsection{The case of biquadratic fields contained in $K$ }
The following example is a  generalization for the biquadratic fields $L$
which are included in $p$-cyclotomic field $K$ with $p\equiv 1\modu 4$.
%
\begin{thm}\label{t607251}
Let $p$ be a prime with $2^2\ \|\ p-1$.
Let
\be\label{e607261}
S=(\sum_{i=0}^{(p-3)/2} (-1)^iv_{2i})^2+(\sum_{i=0}^{(p-3)/2} (-1)^i v_{2i+1})^2.
\ee
Let $L$ be the field with $\Q(\sqrt{p})\subset L\subset K,\ [L:\Q(\sqrt{p})]=2$.
Let $h$ be an odd prime number with $h\ |\ h(L)$ and $h\ \not|\ h(Q(\sqrt{p})$.
Then $S\not=0$ and $S\equiv 0\modu h$.
\begin{proof}
$V(\sigma)\ |\ \sigma^4-1$.
$h\ \not|\ h(\Q(\sqrt{p})$  and so $V(\sigma)\ |\ \sigma^2+1$.
$P(\sigma)=\sum_{i=0}^{(p-3)/2}\sigma^{2i} v_{-2i}+\sigma\times\sum_{i=0}^{(p-3)/2}\sigma^{2i}v_{2i+1}$.
Relation (\ref{e607261}) follows.
\end{proof}
\end{thm}
%
\paragraph{Remarks:}
\bn
\item
$S$ does not depend of the primitive root $v \modu p$ chosen.
\item
Numerical computations seem to show  more : $P(\sigma)\equiv 0\modu p^2$ and so
\be\label{e607262}
\frac{(\sum_{i=0}^{(p-3)/2} (-1)^iv_{2i})^2+(\sum_{i=0}^{(p-3)/2} (-1)^i v_{2i+1})^2}{p^2}\equiv 0\modu h.
\ee
\item
This result is a generalization for biquadratic fields of theorem 145 of Hilbert for quadratic fields.
\en
%
\section{A numerical MAPLE  algorithm}\label{s701081}
This section contains a MAPLE algorithm connected to  the structure of the relative class group $\mathbf C^-$.
For each prime number $p<500$, the algorithm computes
\bn
\item
a primitive root $v\modu p$,
\item
the Stickelberger polynomial $P(X)$  has been defined in lemma \ref{l12161} p. \pageref{l12161} and $Q(X)$
in lemma \ref{l512165} p. \pageref{l512165}.
For the prime numbers $h$ with $3\leq h\leq p^2$ and the primitive root $v\modu p$,
the algorithm  computes the polynomial $GCD(X)\in {\bf F}_h[X]$ for  the  indeterminate $X$
given by formulas:
\bn
\item
if $h\not=p$  then
\item
\be\label{e701051}
GCD(X)=Gcd( P(X), X^{(p-1)/2}+1) \modu h,
\ee
\item
if $h=p$ then
\be
GCD(X)=Gcd(Q(X),X^{(p-1)/2}+1) \modu h,
\ee
\en
\en

The results are compared with the corresponding  tables of  class numbers  of cyclotomic fields $K=Q(\zeta_p)$
for the primes $p<500$ in   Schoof \cite{sch}.
We observe   that:
\bn
\item
The set of odd prime numbers $h$ with $degree(GCD(X))>0$ is  strictly the set of odd prime divisors $h$ of $K$ in Schoof table p. 1142.
\item
The rank $\rho$ of the $h$-Sylow subgroup $\mathbf C(h)$ in Schoof tables 4.2 p. 1239 ($h$ not dividing $p-1$) and
4.3 p. 1240 ($h$ dividing $p-1$) is the degree of the polynomial $GCD(X)$ found here.
Observe that when $h(K^+)\not\equiv 0\modu p$ and  $h=p$ this fact can be proved
by means out of reach of this  article at elementary level (Ribet theorem).
\item
Let $GCD(X)=\prod_{i=1}^n A_i(X)^{n_i}$ be the prime decomposition of $GCD(X)$ in the euclidean field  ${\bf F}_h[X]$.
We observe that to each  prime polynomial $A_i(X)$ corresponds a subgroup of $\mathbf C(h)$ of $h$-rank $d_i\times n_i$ where $d_i$ is the degree of $A_i(X)$.
In particular if $d_i=n_i=1$ then the subgroup corresponding to $A_i(X)$ is cyclic.
Note that for all prime $p$ with $h(K^+)\not\equiv 0\modu p$ and $h=p$ it is a consequence of Ribet theorem.
\item
We observe that when $h\not=p$ there exists some cases where $n_i>1$ for instance for $p=101$ and $h=5$ with $n_1=2$.
In these cases this implies that  the minimal polynomial $V(X)$ annihilating $C_h$   is different of $GCD(X)$ : $V(X)\ |\ GCD(X)$ and
$degree(V(X))< degree (GCD(X))$.
\en
\paragraph{A question:} for all the  odd primes $p<500$ and all the odd primes
$h<p^2$ we have observed  that  the degree of $GCD(X)$ is equal to  the rank $\rho$ of the $h$-Sylow subgroup $\mathbf C(h)$ of the relative class group
$\mathbf C^-$ of $K$
with  $GCD(X)=1\Leftrightarrow  h(K)\not\equiv 0\modu h$: this gives  important informations on the structure of the relative class group
$\mathbf C^-$ of $K$: the precise set of odd primes $h$ dividing $h(K)$ and for each of them the rank $\rho$ of the $h$-group $\mathbf C(h)$.
Can we generalize   this property   to all the odd primes $p$ and all the odd primes $h$ or at least at some  predefined subsets of them?
\paragraph {The MAPLE algorithm}
{\scriptsize
\begin{verbatim}
restart;
> p:=3:
> while p<499 do
> p:=nextprime(p):
> for v from 2 to p-2 do:
> i_v:=1:
> for i from 2 to p-2 do:
> if v&^i mod p = 1 then i_v:=0: fi:
> od:
> if i_v=1 then
> T:= X^((p-1)/2)+1 :
> S:=0:
> Q:=0:
> for i from 0 to p-2 do:
> vmi:=v&^(p-1-i) mod p:
> vmim1:=v&^(p-1-(i-1)) mod p:
> delta_i:=iquo(vmim1-v*vmi,p):
> S:=S+X^i*vmi:
> Q:=Q+X^i*delta_i:
> od:
> for h from 3 to p^2 do:
> if isprime(h)=true then
> if h<>p then
> GCD:=Gcd(T,S) mod h:
> deg_GCD:= degree (GCD):
> if deg_GCD>0  then
> GCD_Factors:=Factors(GCD) mod h:
> n:=nops(GCD_Factors[2]):
> rho:=0:
> for i from 1 to n do:
> temp:=GCD_Factors[2]:
> temp:=temp[i]:
> temp1:=temp[1]:
> temp2:=temp[2]:
> rho:=rho+degree(temp1)*temp2:
> od:
> print(`p=`,p,`h=`,h,`rho=`,rho,`v=`,v,`GCD(X)=`,GCD_Factors[2]):
> fi:
> fi:
> if h = p then
> GCD1:=Gcd(T,Q) mod h:
> deg_GCD1:= degree (GCD1):
> if deg_GCD1>0 then
> GCD1_Factors:=Factors(GCD1) mod h:
> n:=nops(GCD1_Factors[2]):
> rho_1:=0:
> for i from 1 to n do:
> temp:=GCD1_Factors[2]:
> temp:=temp[i]:
> temp1:=temp[1]:
> temp2:=temp[2]:
> rho_1:=rho_1+degree(temp1)*temp2:
> od:
> print(`p=`,p,`h=`,h,`rho_1=`,rho_1,`v=`,v,`GCD1=`,GCD1_Factors[2]):
> fi:
> fi:
> fi:
> od:
> break:
> fi:
> od:
> #fi:
> od:
>
\end{verbatim}}

In the following page the table of results obtained.
\clearpage
{\scriptsize
\begin{tabular}{lllll}
$p$& $h$& $\rho$ &$v$ &$GCD(X)$\\
$23$ & 3 &1&5  & X + 1\\
31&3&1&3&X + 1\\
37&37&1& 2&X + 5\\
41&11&2&6&$X^2  + 10 X + 6$\\
43&211&1&3&X + 73\\
47&5&1&5&X + 1\\
47&139&1&5&X + 44\\
59&3&1&2&X + 1\\
59&59&1&2&X + 36\\
59&233&1&2&X + 8\\
61&41&1&2&X + 36\\
61&1861&1&2&X + 997\\
67&67&1&2&X + 24\\
71&7&1&7&X + 1\\
73&89&1&5&X + 12\\
79&5&1& 3&X + 1\\
79&53&1& 3&X + 28\\
83& 3& 1& 2&X + 1\\
89&113& 1& 3&X + 95\\
97&577& 1& 5&X + 46\\
97& 3457& 1& 5&X + 1558\\
101& 5& 2& 2&$(X + 3)^2$\\
101&101&1& 2&X + 66\\
101& 601 &1& 2&X + 323\\
103& 5& 1& 5&X + 1\\
103& 103& 1& 5&X + 58\\
103& 1021& 1& 5&X + 9\\
107& 3& 1& 2&X + 1\\
107& 743& 1& 2&X + 50\\
107& 9859& 1& 2&X + 4936\\
109& 17& 1& 6&X + 4\\
109& 1009& 1& 6&X + 41\\
113& 17& 1& 3&X + 5\\
127& 5& 1& 3&X + 1\\
127& 13& 1& 3&X + 9\\
127& 43& 1& 3&X + 4\\
127& 547& 1& 3&X + 169\\
127& 883& 1& 3&X + 336\\
127& 3079& 1& 3&X + 1925\\
131& 3& 3& 2&$X^3  + 2 X^2  + 1$\\
131& 5&1& 2&X + 1\\
131& 53& 1& 2&X + 46\\
131& 131&1& 2&X + 34\\
131&1301& 1& 2&X + 283\\
137&17&1&3&X + 8\\
139& 3& 1& 2&X + 1\\
139&47& 1& 2&X + 9\\
139& 277& 2& 2&(X + 191)(X + 218)\\
139& 967& 1& 2&X + 241\\
149& 3& 2& 2&$X^2  + 1$\\
149& 149& 1& 2&X + 43\\
151& 7& 1& 6&X + 1\\
151& 11& 2& 6& $X^2 + 6 X + 3$\\
151&281& 1& 6&X + 90\\
\end{tabular}
\begin{tabular}{lllll}
$p$& $h$& $\rho$ &$v$ &$GCD(X)$\\
157& 5&1& 5&X + 3\\
157& 13& 1&5&X + 6\\
157& 157& 2& 5&(X + 95)(X + 91)\\
157& 1093& 1& 5&X + 800\\
157&1873& 1& 5&X + 935\\
163&181& 1& 2&X + 65\\
163& 23167&1& 2&X + 8783\\
167& 11& 1& 5&X + 1\\
167& 499& 1& 5&X + 491\\
173& 5&1&2&X + 2\\
173& 20297& 1& 2&X + 997\\
179& 5& 1& 2&X + 1\\
179& 1069& 1& 2&X + 552\\
181& 5& 1& 2&X + 3\\
181& 37& 1& 2&X + 29\\
181& 41& 1& 2&X + 2\\
181& 61& 1& 2&X + 6\\
181& 1321& 1& 2&X + 149\\
181&2521& 1& 2&X + 2015\\
191& 11& 1& 19&X + 3\\
191& 13& 1& 19&X + 1\\
193& 6529& 1& 5&X + 4193\\
193&15361& 1& 5&X + 13057\\
193& 29761& 1& 5&X + 29163\\
197& 5& 1& 2&X + 3\\
197& 1877& 1& 2&X + 981\\
197& 7841& 1& 2&X + 1604\\
199& 3& 1& 3&X + 1\\
199& 19& 1& 3&X + 4\\
199& 727& 1&3&X + 590\\
211& 3& 1& 2&X + 1\\
211& 7& 1& 2&X + 4\\
211& 41& 1& 2&X + 16\\
211& 71& 1& 2&X + 15\\
211& 181& 1& 2&X + 5\\
211& 281& 2& 2&(X + 109)(X + 101)\\
211& 421&1& 2&X + 93\\
211& 1051& 1& 2&X + 884\\
211& 12251& 1& 2&X + 1580\\
223& 7&1& 3&X + 1\\
223& 43& 1& 3&X + 36\\
227& 5& 1&2&X + 1\\
227& 2939& 3& 2&(X + 1420)(X + 509)(X + 2006)\\
229&13&1& 6&X + 6\\
229& 17& 1& 6&X + 13\\
229& 457& 1& 6&X + 126\\
229& 7753& 1& 6&X + 4310\\
233& 233& 1& 3&X + 193\\
233& 1433& 1& 3&X + 1091\\
239& 3& 1& 7&X + 1\\
239& 5& 1& 7&X + 1\\
\end{tabular}
\clearpage
\begin{tabular}{lllll}
$p$& $h$& $\rho$ &$v$ &$GCD(X)$\\
241& 47& 2& 7&$X^2  + 29 X + 1$\\
241& 13921& 1& 7&X + 9953\\
241& 15601& 1& 7&X + 7049\\
251& 7& 1& 6&X + 1\\
251& 11&  1& 6&X + 4\\
257& 257& 1& 3&X + 76\\
263& 13& 1& 5&X + 1\\
263& 263& 1& 5&X + 204\\
263& 787& 1& 5&X + 510\\
269& 13& 1& 2&X + 5\\
271& 11& 1& 6&X + 1\\
271& 31& 1& 6&X + 2\\
271& 37& 1& 6&X + 9\\
271& 271&1& 6&X + 196\\
271& 811& 1& 6&X + 133\\
271& 1201& 1& 6 &X + 367\\
271& 1621& 1& 6&X + 1190\\
271& 15391& 1& 6&X + 6331\\
271& 21961& 1& 6&X + 6698\\
277& 17& 1& 5&X + 13\\
277& 47&  2& 5&$X^2  + 19 X + 42$\\
277& 829& 1& 5&X + 150\\
281& 11& 2& 3&$X^2  + 2 X + 2$\\
281& 17& 1& 3&X + 8\\
281& 41& 1& 3&X + 24\\
281& 401&1& 3&X + 250\\
283& 3& 1& 3&X + 1\\
283& 283& 1& 3&X + 236\\
293& 3& 2& 2&$X^2  + 1$\\
293& 293&1& 2&X + 194\\
307& 3& 1& 5&X + 1\\
307& 37& 1& 5&X + 33\\
307& 137& 1& 5&X + 38\\
307& 307& 1& 5&X + 16\\
307& 443& 1& 5&X + 13\\
307& 613& 1& 5&X + 49\\
307& 919& 1& 5&X + 144\\
307& 1429& 1& 5&X + 1294\\
311& 19& 1& 17&X + 1\\
311& 41& 1& 17&X + 18\\
311& 311& 1&17&X + 158\\
313& 37& 2& 10&$X^2  + 14$\\
313& 233& 1& 10&X + 136\\
317& 13& 1& 2&X + 5\\
331& 3& 6& 3&$(X + 1)^2(X^4  + 2 X^3  + X ^2 + 2 X + 1)$\\
331& 23& 1& 3&X + 2\\
331& 61& 1& 3&X + 12\\
331& 67&1& 3&X + 64\\
337& 7& 2& 10&$X^2  + X + 6$\\
337& 17& 2& 10&(X + 10)(X + 7)\\
337& 353& 1& 10&X + 36\\
347& 5& 1& 2&X + 1\\
347&347& 1& 2&X + 52\\
349& 5& 1& 2&X + 2\\
349& 13& 1& 2&X + 6\\
349& 2089& 1& 2&X + 1733\\
349& 17749& 1& 2&X + 9289\\
\end{tabular}
\begin{tabular}{lllll}
$p$& $h$& $\rho$ &$v$ &$GCD(X)$\\
353& 353& 2& 3&(X + 299)(X + 51)\\
353& 6113& 1& 3&X + 2060\\
353& 9473& 1& 3&X + 5067\\
359& 19& 1&7&X + 1\\
367& 3& 1& 6&X + 1\\
367& 733& 1& 6&X + 686\\
367& 39163& 1& 6&X + 27454\\
373& 5& 1& 2&X + 2\\
373& 61& 1& 2&X + 21\\
373&1117&1& 2&X + 532\\
373& 1489& 1& 2&X + 990\\
379& 3& 1& 2&X + 1\\
379& 13& 1& 2&X + 3\\
379& 127& 1& 2&X + 13\\
379& 379& 2& 2&(X + 348)(X + 91)\\
379& 547& 1& 2&X + 196\\
379& 757& 1& 2&X + 531\\
379& 991& 1& 2&X + 324\\
379& 1499& 1& 2&X + 1314\\
379& 9199& 1& 2&X + 7901\\
383& 17& 1& 5&X + 1\\
389& 41& 1& 2&X + 32\\
389& 389& 1& 2&X + 231\\
389& 1553& 1& 2&X + 1130\\
397& 13&1& 5&X + 5\\
397& 23& 2& 5&$X^2  + 13 X + 18$\\
397& 109& 1& 5&X + 32\\
397&4861& 1& 5&X + 3655\\
397& 9901& 1& 5&X + 8544\\
401&41& 2& 3&$X^2  + 17$\\
401& 401& 1& 3&X + 141\\
401& 64849& 1& 3&X + 46775\\
409& 5& 2& 21&$X^2  + 3$\\
409& 17& 1& 21&X + 2\\
409& 73& 1& 21&X + 66\\
409& 409& 1& 21&X + 28\\
409& 1321& 1& 21&X + 1304\\
419& 3& 1& 2&X + 1\\
419& 1103& 1& 2&X + 494\\
421& 5& 1& 2&X + 2\\
421& 29& 1& 2&X + 18\\
421& 37& 1& 2&X + 23\\
421& 421& 1& 2&X + 72\\
421& 2521&1& 2&X + 60\\
421& 39509& 1& 2&X + 7582\\
421& 39901& 1& 2&X + 8081\\
421& 70309& 1& 2&X + 65038\\
431& 3& 1& 7&X + 1\\
431& 7&1&7&X + 1\\
431& 11& 1& 7&X + 5\\
431& 701& 1& 7&X + 210\\
431& 14621& 1& 7&X + 2522\\
433& 433& 1& 5&X + 371\\
433& 3457& 1& 5&X + 2700\\
433& 12097& 1& 5&X + 31\\
433& 21601& 1& 5&X + 10658\\
433& 47521& 1& 5&X + 36247\\
\end{tabular}
\clearpage
\begin{tabular}{lllll}
$p$& $h$& $\rho$ &$v$ &$GCD(X)$\\
439& 3& 2& 15&$(X + 1)^2$\\
439& 5& 1&15&X + 1\\
439& 293& 1& 15&X + 283\\
443& 3& 3& 2&$X^3  + X^2  + 2 X + 1$\\
443& 5&1& 2&X + 1\\
443& 79& 1& 2&X + 8\\
443& 157& 1& 2&X + 67\\
443& 12377& 1& 2&X + 6026\\
449& 168449& 1& 3&X + 33570\\
457& 5& 2& 13&$X^2  + 4 X + 2$\\
457&41& 1& 13&X + 14\\
457& 577& 1& 13&X + 9\\
457& 1217& 1& 13&X + 692\\
457& 43777& 1& 13&X + 37577\\
457& 63841& 1& 13&X + 2827\\
461& 5& 2&2&$(X + 2)^2$\\
461& 461&1& 2&X + 13\\
461& 661& 1& 2&X + 258\\
461& 161461& 1&2&X + 134936\\
463& 7& 2& 3&(X + 2)(X + 1)\\
463& 29& 1& 3&X + 20\\
463& 89& 1& 3&X + 64\\
463& 463&1& 3&X + 8\\
463& 631& 1& 3&X + 62\\
463& 673& 1&3&X + 223\\
463& 1123& 1& 3&X + 49\\
463& 4423& 1& 3&X + 387\\
463& 8779& 1& 3&X + 5520\\
467& 7& 1& 2&X + 1\\
467& 467& 2& 2&(X + 239)(X + 236)\\
479& 5& 1&13&X + 1\\
479& 48757& 1& 13&X+34844\\
479& 62141& 1& 13&X + 43049\\
487& 7& 2& 3&(X + 4)(X + 1)\\
487& 37& 2& 3&(X + 33)(X + 12)\\
487& 919& 1& 3&X + 267\\
487& 2647& 1& 3&X + 1070\\
487& 10909& 1& 3&X + 3031\\
487& 58321& 1& 3&X + 58241\\
491& 3& 1& 2&X + 1\\
491& 11& 2& 2&(X + 4)(X + 5)\\
491& 29& 1& 2&X + 16\\
491& 491& 3& 2&(X + 203)(X + 419)(X + 418)\\
499& 3& 1& 7&X + 1\\
499& 167& 1& 7&X + 98\\
\end{tabular}}

%

%
Roland Qu\^eme

13 avenue du ch\^ateau d'eau

31490 Brax

France

mailto: roland.queme@wanadoo.fr

home page: http://roland.queme.free.fr/

************************************

V10 - MSC Classification : 11R18;  11R29

************************************
\end{document}